\documentclass[12pt]{amsart}
\usepackage{a4,amsmath,amssymb,amsthm,eucal}
\usepackage[english]{babel}
\usepackage{graphicx}
\usepackage{xcolor}
\usepackage{epsfig}

\newcommand{\forget}[1]{}

\def\Om{{\Omega}}  
\def\om2{{\Om\times\Om}}

\usepackage{color}

\newtheorem{theorem}{Theorem}[section]
\newtheorem{definition}[theorem]{Definition}

\newtheorem{lemma}[theorem]{Lemma}
\newtheorem{proposition}[theorem]{Proposition}

\theoremstyle{remark}
\newtheorem{remark}[theorem]{Remark}

\numberwithin{equation}{section}

\usepackage[normalem]{ulem}

\definecolor{blue-violet}{rgb}{0.54,0.17,0.89}
\definecolor{amethyst}{rgb}{0.6,0.4,0.8}
\definecolor{darkviolet}{rgb}{0.58, 0.0, 0.83}
\definecolor{darkgreen}{rgb}{0,.4,0}
\definecolor{mixedgreen}{rgb}{0.3,0.6,00}
\definecolor{bananayellow}{rgb}{1.0, 0.88, 0.21}
\definecolor{arylideyellow}{rgb}{0.91, 0.84, 0.42}
\definecolor{bananamania}{rgb}{0.98, 0.91, 0.71}

\newcommand{\REMOVEsk}[1]%
           {{\color{magenta}\sout{#1}}}

\begin{document}
\title[Bounded and exponentially decaying solutions]{Bounded and exponentially decaying solutions of almost linear dynamic systems on time scales.}
\author{Svetlin Georgiev}
\author{Sergey Kryzhevich}
\address[Sergey Kryzhevich]{
Institute of Applied Mathematics, Faculty of Applied Physics and Mathematics 
\and
BioTechMed Center, Gda\'nsk University of Technology, 80-233 Gda\'nsk, Poland}
\email[Sergey Kryzhevich]{serkryzh@pg.edu.pl}
\address[Svetlin Georgiev]{Sorbonne University, 19 Rue de l'Ecole de M\'edecine, 75006 Paris, France}
\email[Svetlin Georgiev]{svetlingeorgiev1@gmail.com}

\thanks{The work was supported by Gda\'{n}sk University of Technology by the DEC 14/2021/IDUB/I.1 grant under the Nobelium - 'Excellence Initiative - Research University' program.}


\begin{abstract}
We consider small nonlinear perturbations of linear systems on a time scale with the phase space being finite or infinite-dimensional. For $\Delta$-differential operators, corresponding to linear dynamic systems we consider their solvability in various functional spaces. Based on these techniques, we prove several results on the existence of a bounded solution for some systems with small nonlinearities. Thus we introduce some generalizations of the classic hyperbolicity property (exponential dichotomy) for systems of ordinary differential equations. 

Besides, we introduce the Lyapunov regularity condition that coincides with the classic one for ordinary differential equations. For regular systems, we prove some criteria for the existence of bounded solutions of perturbed systems.  
\end{abstract}

\maketitle

\noindent{\small {\bf Keywords:} time scale dynamics; structural stability; solvability; linear systems; small nonlinear perturbations; Perron approximations, Hilger exponents}

\section{Introduction}

Time scale dynamics is a very popular topic in the contemporary theory of dynamical systems. First introduced by Aulbach and Hilger \cite{AH90} these systems have been intensively studied during the last decades. The survey of the results obtained in the area can be found in \cite{BP03}, \cite{G19}, and \cite{M16}. Undoubtedly, the study of various classes of perturbations is of paramount importance for time scale systems as well as those for ordinary differential equations or discrete dynamical systems.

Concerning the individual stability of solutions, several analogs of classical approaches have been developed based on the Gronwall -- Bellman inequality, Lyapunov functions theory, and hyperbolicity theory. However, there appear significant difficulties in translating the classic results of structural stability and shadowing theory to the language of time scales. 

First of all, generic time scales are non-periodic and, therefore, one can hardly expect the existence of periodic solutions (and even if they exist, they do not play any specific role in the dynamics of the system). This is why it is difficult to formulate analogs of classic conditions (for example, those of Axiom A) for generic time scale systems, similar to non-autonomous ordinary differential equations.

Secondly, even the uniqueness of solutions is a non-trivial property for time scale systems. In other words, it does not follow directly from the smoothness of the right-hand side, some additional assumptions have to be imposed. 

Finally, there is a more specific feature of time scale dynamics. The time scales with unbounded graininess function (roughly speaking, those with infinitely large gaps in the scale) may admit no linear systems with exponentially increasing/decreasing non-trivial solutions. This may lead to difficulties in defining hyperbolicity for such systems. That is why many authors, dealing with time scale systems prefer to work with the so-called Hilger exponents (see Section 4) rather than regular ones.    

In this paper, we suggest the following approach. Given a space of right-hand sides of linear nonhomogeneous systems with a fixed homogeneous part, we find conditions on the homogeneous part sufficient for the considered space of systems to have a bounded solution. For instance, it is well-known that a hyperbolic hon-homogeneous system of ODEs with a bounded right-hand side admits a bounded solution while a regular system with an exponentially decaying right-hand side has an exponentially decaying solution. In this paper, we study some general classes of time scale systems on Banach phase space. In addition, we generalize the concept of hyperbolicity even without assuming that the dynamics are invertible. For example, we may take a non-regressive linear system with $\Delta$ -- derivative. This may imply the absence of a backward uniqueness of solutions which is treated as an infinite shrink of solutions of a certain subspace of the phase space. 

The rest of the paper consists of the following parts. 

In Section 2, some basic facts of time scale calculus and dynamics are recalled.

The most general situation is considered in Section 3 where some classes of 'weak' hyperbolicity are introduced and some related results are proved. Besides, we consider the specific case of finite-dimensional phase spaces. Under this assumption, the condition of the smallness of the Lipschitz coefficient for the right-hand side can be slightly weakened.

In Section 4, we introduce the theory of Lyapunov exponents and regularity via Hilger exponents of a specially selected family of functions. A brief theory of these exponents and their stability is also given.

In Section 5 we spread the concept of Lyapunov regular systems to the case of time scale dynamics and study how this kind of regularity is connected with solvability for linear nonhomogeneous systems. 

The conclusion and discussion are given in Section 6.

\section{Some basic definitions and facts from Time Scale Dynamics}

In this section, we use the notions of books \cite{BP03} and \cite{G19}. 

The main object of our interest is a closed non-empty subset ${\mathbb T}\subset {\mathbb R}$ with the inherited topology. We call any of such sets \emph{time scale}. 

\begin{definition}
Given a $t \in {\mathbb T}$ one defines the so-called \emph{forward jump operator} $\sigma : {\mathbb T} \mapsto {\mathbb T}$ by the formula
$\sigma(t) = \inf\{s \in {\mathbb T} : s > t\}$. Observe that $\sigma(t) \ge t$ for any $t \in {\mathbb T}.$ If the above set is empty, we set $\sigma(t)=\sup {\mathbb T}$.
\end{definition}

\begin{definition}\label{rdrs}
We define the graininess function as follows $\mu(t)=\sigma(t)-t$. The point $t\in {\mathbb T}$ is called \emph{right-dense} if $\mu(t)=0$ and right-scattered otherwise. Backward jump operators, left-dense and left-scattered points can be defined similarly.
\end{definition}

For every time scale $\mathbb T$ the set $RS({\mathbb T})$ of all right-scattering points is countable.

\begin{definition}
We say that a function is rd-continuous if it is continuous in all right-dense points and there exists a left limit at all left-dense points $($those limits do not necessarily coincide with the values of the function$)$.
\end{definition}

\begin{definition} If $f : {\mathbb T} \mapsto {\mathbb R}$ is a function, then we define the function
$f^\sigma : {\mathbb T} \mapsto {\mathbb R}$ by
$f^{\sigma}(t)=f(\sigma(t))$ for any $t \in {\mathbb T}$, i.e., $f^\sigma =f\circ \sigma$.
\end{definition}

\begin{definition}
Assume that $f : {\mathbb T} \mapsto {\mathbb R}$ is an rd-continuous function. Given a $t \in {\mathbb T}$ we define the so-called Hilger $\Delta$-derivative (or just $\Delta$-derivative) as follows. We say that $f^\Delta(t)=A$ if for any $\varepsilon > 0$ there is a neighborhood $U$ of $t$ in $\mathbb T$, such that
$$|f^{\sigma}(t)-f(s)-A (\sigma(t)-s)|\le \varepsilon|\sigma(t)-s|\quad \mbox{for all}\quad s\in U, s \neq \sigma(t).$$

We say that $f$ is differentiable, in $\mathbb T$ if $f^\Delta (t)$ exists for all $t \in {\mathbb T}$. The function $f^\Delta : {\mathbb T} \mapsto {\mathbb R}$ is called the derivative of $f$ in this case.
\end{definition}

For ${\mathbb T}={\mathbb R}$, the derivative is the classical one, if ${\mathbb T} = {\mathbb N}$, we have $f^\Delta(n)=f(n+1)-f(n)$.

We list some basic properties of Hilger derivatives.

\begin{theorem} Let $f :{\mathbb T} \mapsto R$ be a function and let $t\in {\mathbb T}$. 
\begin{enumerate}
\item If $f$ is differentiable at $t$, then $f$ is continuous at $t$.
\item If $f$ is continuous at $t$ and $t$ is right-scattered, then $f$ is differentiable at $t$ with
$$f^\Delta(t)= \dfrac{f^{\sigma}(t)-f(t)}{\mu(t)}.$$
\item If $t$ is right-dense, then $f$ is differentiable if and only if the limit
$$f^\Delta(t):=\lim_{s\to t} \dfrac{f(t)-f(s)}{t-s}$$
exists and is finite.
\item If $f$ is differentiable at $t$, then $f (\sigma (t)) = f (t) + \mu(t)f^\Delta (t)$.
\end{enumerate}
\end{theorem}

\begin{theorem} For any functions $f,g:{\mathbb T}\mapsto {\mathbb R}$ differentiable at $t\in {\mathbb T}$ 
\begin{enumerate}
\item the sum $f+g$ is differentiable at $t$ and
$$(f+g)^\Delta(t)=f^\Delta(t)+g^\Delta(t).$$
\item for any constant $\alpha$, the function $\alpha f$ is differentiable at $t$ and $$(\alpha f )^\Delta(t) = \alpha f^\Delta(t).$$
\item the product $fg$ is differentiable at $t$ and
$$(fg)^\Delta(t) = f^\Delta(t)g(t)+f^{\sigma}(t)g^\Delta(t)
= f(t)g^\Delta(t)+f^\Delta(t)g^{\sigma}(t).$$
\item if $g(t),g^{\sigma}(t)\neq 0$ then $f/g$ is differentiable at $t$ and
$$(f/g)^\Delta (t)= \dfrac{f^\Delta(t)g(t)-f(t)g^\Delta(t)}{g(t)g^{\sigma}(t)}.$$
\end{enumerate}
\end{theorem}
The antiderivatives can be defined on time scales similar to classic calculus. Then, the definite integrals
$$\int_a^b f(t)\, \Delta t$$
can be introduced by the Newton-Leibnitz formula which is used as a definition.

One can consider analogs of differential equations and, in particular, linear systems of the type
\begin{equation}\label{linear}
    x^\Delta=A(t)x
\end{equation}
on time scales. In general, the properties of systems \eqref{linear} are very similar to those for ordinary differential equations with one important exception which can be illustrated by the following example.

\noindent\textbf{Example}. Let ${\mathbb T}={\mathbb N}$. Consider the scalar equation
$x^\Delta=-x$. Then any solution of the considered equation with initial conditions $x(n_0)=x_0$ is zero for any $n>n_0$, so there is no backward uniqueness of solutions. Moreover, if $x_0\neq 0$, the solutions of the considered Cauchy problem are not defined for $n<n_0$.

\begin{definition} \hfill
\begin{enumerate}
\item We say that an $rd$-continuous function $a:{\mathbb T}\to {\mathbb C}$ is regressive if $1+\mu(t)a(t)\neq 0$ for all $t\in {\mathbb T}$.
\item This function is positively regressive if $1+\mu(t)a(t)>0$ for all $t\in {\mathbb T}$.
\item This function is uniformly positively regressive if $\inf_{\mathbb T}(1+\mu(t)a(t))>0$.
\item We say that the matrix $A$ is regressive with respect to $\mathbb T$ provided $E+\mu(t)A(t)$ is invertible for all $t\in {\mathbb T}$.
\end{enumerate}
\end{definition}
Denote the classes of regressive, positively regressive, and uniformly regressive functions by $\mathcal R$, ${\mathcal R}^+$ and 
${\mathcal{UR}}^+$ respectively.

\begin{theorem}\hfill
\begin{enumerate}
    \item Any solution of system \eqref{linear} with initial conditions $x(t_0)=x_0$ exists and is uniquely defined for any $t\ge t_0$.
    \item Solutions of system \eqref{linear} are unique if the matrix $A(t)$ is regressive.
\item The matrix-valued function $A$ is regressive if and only if the eigenvalues $\lambda_i(t)$ of $A(t)$ are regressive for all $1 \le  i \le  n$.
\end{enumerate}
\end{theorem}

An analog of the Lagrange method to solve linear non-homogeneous systems holds for time scale systems, as well. Namely, given a system \eqref{linear} on a Banach space $X$ and
$f:{\mathbb L}^1_{loc} ({\mathbb T}\mapsto X)$, we consider the non-homogeneous system
\begin{equation}\label{timescale3}
x^\Delta=A(t)x+f(t).
\end{equation}
with a bounded matrix $A(t)$ (if $X$ is infinite-dimensional we still call the operator $A(t)$ matrix). Let $\Phi(t,s)$ be the Cauchy matrix for system \eqref{linear}, i.e., the fundamental matrix of \eqref{linear} with $\Phi(s,s)=E$. Here $E$ stands for the unit matrix/operator. Generally speaking, such a matrix is defined for any $t\ge s$. Then the solution of problem \eqref{timescale3} with initial conditions $x(t_0)=x_0$ may be found by the following formula
$$x(t)=\Phi(t,t_0)x_0+\int_{t_0}^t \Phi(t,\sigma(s))\, \Delta s.$$
Observe that this formula is true even for non-regressive systems.

Let us formulate the statement we will need later on.

\begin{lemma}\label{non-unif} Let 
$t_-=\inf {\mathbb T}$, $t_+=\sup {\mathbb T}$ where both the values may be finite or infinite.
Consider two $rd$-continuous and uniformly bounded families of projections $P(t)$ and $Q(t)$ such that $P(t)+Q(t)=E$.
Then
$$\int_{t_-}^t \Phi(t,\sigma(s)) P(s)f(s)\, \Delta s -\int_{t}^{t_+} \Phi(t,\sigma(s)) Q(s)f(s)\, \Delta s$$
is a solution of \eqref{timescale3}
\end{lemma}

The proof of the lemma is straightforward.

\section{General case. Linear systems on time scales and their solvability.}

Now we proceed to the case of a generic time scale (a closed non-empty subset of the real axis), still denoted as $\mathbb T$.  

We consider a non-linear system
\begin{equation}\label{timescale1}
x^\Delta=A(t) x + a(t,x)
\end{equation}
where $A(t)$ and $a(t,x)$ are defined and rd-continuous on ${\mathbb T}$ and ${\mathbb T}\times X$, respectively.  Also, we consider the 'linear part' \eqref{linear} of system \eqref{timescale1}.

We suppose that the following condition is satisfied.

\medskip

\noindent\textbf{Condition I.}
\begin{enumerate}
\item The matrix $A$ is bounded.
\item Solutions of system \eqref{timescale1} with any initial conditions $(t_0,x_0)\in {\mathbb T}\times X$ are unique over $\mathbb T\bigcup [t_0,\infty)$.
\end{enumerate}
\medskip

Let $\Phi(t,s)$ be the Cauchy operator for system \eqref{linear}, i.e. the fundamental matrix with $\Phi(s,s)=E$.

Consider a countable family of Banach spaces $(V_j,\|\cdot\|_{V_j})$, $j\in J$, and a Banach space $(W,\|\cdot\|_W)$ all consisting of functions from the time scale $\mathbb T$ to the space $X$.

\medskip

\noindent\textbf{Condition II.}
\begin{enumerate}
    \item There exists a continuous family of bounded projections $P_j(t)$ and $Q_j(t)$ of the space $X$ ($j\in J$, $t\in {\mathbb T}$) such that $P_j(t)+Q_j(t)=E$ for any $j$ where $E$ stands for the identity operator on the space $X$ and, given a right-hand side $f$, the function ${\mathcal L}_j f(t)$ defined by the formula
    \begin{equation}\label{eq_lj}
    {\mathcal L}_j f(t)=
    \int_{t_-}^t \Phi(t,\sigma(s)) P_j(s)f(s)\, \Delta s -\int_{t}^{t_+} \Phi(t,\sigma(s)) Q_j(s)f(s)\, \Delta s
    \end{equation}
    gives a solution of system \eqref{timescale3}. 
    \item Linear operators ${\mathcal L}_j$ defined by Eq.\eqref{eq_lj} act continuously from each of $V_j$ to $W$. We introduce $L_j:=\|{\mathcal L}_j\|$.
    \item There exist positive constants $h_j, c_j: j\in J$ such that 
    $$a(t,x)=\sum_{j\in J} g_j(t,x)$$
    with $g_j (\cdot,x(\cdot)) \in V_j$ for any function $x(\cdot)\in W$ and 
    \begin{equation}\label{nonlinest}
    \|g_j(\cdot,0)\|_{V_j}\le h_j, \quad \|g_j(\cdot,x(\cdot))-g_j(\cdot,y(\cdot))\|_{V_j}\le c_j \|x(\cdot)-y(\cdot)\|_W.
    \end{equation}
    \item The values 
    \begin{equation}\label{elambda}
    \beta:=\sum_{j\in J} L_jh_j <+\infty, \qquad \lambda=\sum_{j\in J} L_jc_j<1.
    \end{equation}
\end{enumerate}

\medskip

In this case, we call the linear system \eqref{linear} $(\{V_j\}, W)$ -- \emph{hyperbolic} and the perturbation $a(t,x)$ $(\{V_j\}, W)$ -- \emph{small}.

\begin{definition}\label{defhyp}
We say that the system is \emph{hyperbolic} on all the time scale if Condition 2 is satisfied with the set $\{V_j\}$ of Banach spaces consisting of one space ${\mathbb L}^\infty({\mathbb T}\to X)$ and with the space $W$ being the same.
\end{definition}

\begin{theorem}\label{tmain} Let system \eqref{timescale1} satisfy Conditions I and II. Then, this system admits a solution $x^*(\cdot)\in W$ with 
\begin{equation}\label{est_sol}
    \|x^*(\cdot)\|_W \le \dfrac{\beta}{1-\lambda}.
\end{equation}
\end{theorem}

\noindent\textbf{Proof.} We apply the Banach fixed point theorem. Indeed, we can define an operator $T:W\mapsto W$ by the formula 
\begin{equation}\label{eqT}
Tf=\sum_{j\in J} {\mathcal L}_j a_j(\cdot,f).
\end{equation}
Inequalities \eqref{elambda} imply that $$\|Tx\|\le \lambda \|x\|+\beta$$ for any $x\in W$ and the operator $T$ is contracting which implies Eq.\,\eqref{est_sol}. $\square$

\medskip

\noindent\textbf{Remark.} We can suppose that inequalities \eqref{nonlinest} hold true for $x:$ $\|x\|\le \beta/(1-\lambda)$ only. Of course, an analog of the proof of the above theorem is still valid in these assumptions.

\bigskip

\noindent\textbf{The finite-dimensional phase space $X$ and hyperbolicity}.
Let $\dim X<+\infty$ and system \eqref{linear} be hyperbolic (exponentially dichotomic). In this case, we assume that the nonlinearity satisfies the following condition instead of \eqref{nonlinest}:
\begin{equation}\label{alinfty}
\|a(t,x)\|\le \lambda_0 \|x\|+\beta_0
\end{equation}

\begin{theorem}\label{th_hyp_1} Let system \eqref{linear} be hyperbolic (see Definition \ref{defhyp}) and the nonlinearity $a(t,x)$ satisfies condition \eqref{alinfty} with a sufficiently small $\lambda_0$. Then system \eqref{timescale1} has a bounded solution. 
\end{theorem}

\noindent\textbf{Proof.} We mainly repeat the proof of Theorem \ref{tmain}. In the considered case we have only one operator $\mathcal L$ acting from ${\mathbb L}^\infty ({\mathbb T})$ to itself. 

Let $\lambda=\lambda_0 \|{\mathcal L}\|$, $\beta=\beta_0 \|{\mathcal L}\|$. Suppose that $\lambda_0$ is so small that $\lambda<1$. The operator $T$, defined by formula
$$Tf= {\mathcal L} a(\cdot,f),$$
similar to \eqref{eqT}, transfers the ball $B$ of continuous functions with $\|f\|\le \beta/(1-\lambda)$ to its subset. Moreover, all the functions of the set $TB$ are uniformly continuous (their derivatives are uniformly bounded). Then, due to Arzel\`a -- Ascoli Theorem, the image $TB$ is a compact subset of $B$ in the topology of uniform convergence on compact subsets of ${\mathbb T}$. Then the Schauder fixed point theorem implies the existence of a bounded solution. $\square$

Observe that we need our phase space to be finite-dimensional in order to be able to apply Arzel\`a -- Ascoli Theorem.

The last statement may be generalized.

\begin{theorem}\label{th_hyp_2} Let $t_-=-\infty$, $t_+=\infty$. Suppose that there exists a $t_0\in {\mathbb T}$ such that system \eqref{linear} is hyperbolic $($also in the sense of Definition \ref{defhyp}$)$ on both the sets 
$${\mathbb T}_-:=(-\infty, t_0] \bigcap {\mathbb T},\quad \mbox{and} \quad {\mathbb T}_+:=[t_0,\infty) \bigcap {\mathbb T}.$$
Let $P_\pm(t)$ and $Q_\pm(t)$ be the corresponding projections. Suppose that the spaces 
$P_-(t_0)X$, and $Q_+(t_0)X$ intersect transversely. 
Assume that the nonlinearity $a(t,x)$ satisfies condition \eqref{alinfty} with a sufficiently small $c$. Then system \eqref{timescale1} has a bounded solution. 
\end{theorem}

Proofs of the two latter results are very similar to those of ordinary differential equations \cite{C78}, \cite{KV06}.

\section{\bf Hilger functions and Lyapunov exponents}

\begin{definition} We say that the time scale $\mathbb T$ is \emph{syndetic} if 
$$
\mu^*=\sup_{t\in {\mathbb T}} \mu(t)<+\infty.
$$
\end{definition}

Starting from this point, we always assume that the considered time scale $\mathbb T$ is unbounded ($\sup {\mathbb T}=+\infty$) and syndetic. Let $\nu^*=1/\mu^*$ if $\sup ({\mathbb R}\setminus {\mathbb T})>0$ and $\nu^*=\infty$ otherwise. Besides, we assume that the phase space $X$ is finite-dimensional, $\dim X=n$.

Consider the solutions $e_p(t,s)$ of Cauchy problems  
\begin{equation}\label{equp}
x^\Delta=p(t) x, \qquad x(s)=1
\end{equation}
where $p$ is an $rd$--continuous scalar function. These functions play an important role in the theory of time scale systems, see \cite{AH90}, \cite{BP01}, and \cite{BP03} for details.
The function $e_p(t,s)$, also known as Hilger exponent, could be found as follows.

Consider the so-called cylinder transformation 
$$\xi_h(z)=\dfrac1{h} \log (1+zh) \qquad \mbox{if} \quad h\neq 0$$
and $\xi_0(z)=z$.
Here $\log$ is the principle logarithm function defined for 
$$-\dfrac{\pi}{h}<\mbox{Im}\, z \le \dfrac{\pi}{h}.$$
If $zh>-1$, $\log$ stands for classic real logarithm. In some papers, e.g. \cite{AB21}, the multi-valued logarithms are considered, and, respectively, the multivalued exponent functions are introduced. Besides, alternative definitions of exponents on time scales are given in \cite{C12}.   

Then 
$$e_p(t,s)=\exp\left(\int_s^t\xi_{\mu(\tau)}\, p(\tau)\, \Delta \tau\right).$$
unless there exists a $t_0\ge s$ such that $\mu(t_0)p(t_0)=-1$. In that case $e_p(t,s)=0$ for any $t\ge t_0$. If ${\mathbb T}={\mathbb R}$, 
$$e_p(t,s)=\exp\left(\int_{s}^t p(\tau)\, d \tau\right).$$

Anyway, the function $e_p(t,s)$ is positive if $p(t)$ is positively regressive. In particular, if $p=\mbox{const}$, the function $e_p(t,s)$ is positive if $p>-1/\mu^*$. By contrast, if $p<-1/\mu^*$, then the function $e_p(t,s)$ is non-positive for some values of $t$ and $s$. 

These exponents are closely connected with the following operations on the set of functions (these definitions are also classic):
$$\begin{array}{c}
(p\oplus q)(t)=p(t)+q(t)+\mu (t) p(t)q(t), \quad (p\ominus q)(t)=\dfrac{p(t)-q(t)}{1+\mu (t) q(t)}, \\
\ominus q (t)=\dfrac{-q(t)}{1+\mu (t) q(t)}.
\end{array}$$
These operators define a commutative group structure on the set of functions (but not on the set of real/complex numbers). Observe that given a function $p(t)\in {\mathcal{UR}}^+$ and any $\varepsilon>0$ there exists an $\varepsilon_1>0$ such that
$(p\oplus \varepsilon)(t) \ge p(t) + \varepsilon_1$

Let us list some simple and well-known properties of Hilger exponents (see \cite{BP03}).

\begin{proposition} Given two $rd$ -- continuous functions $p(\cdot)$ and $q(\cdot)$, and any $t,s,\tau \in {\mathbb T}$, we have  
\begin{enumerate}
\item $e_p(t,t)=0$, $e_0(t,s)=1$.
\item $e_p(t,s)e_p(s,\tau)=e_p(t,\tau)$.
\item $e_p(\sigma(t),s)=(1+\mu(t)p(t))e_p(t,s)$.
\item $e_p(t,s)e_q(t,s)=e_{p\oplus q}(t,s).$
\item $e_p(t,s)/e_q(t,s)=e_{p\ominus q}(t,s),$ $1/e_q(t,s)=e_q(s,t)=e_{\ominus q}(t,s)$.
\item If $p(\cdot), q(\cdot) \in {\mathcal R}^+$ and $p\le q$ then 
$$0<e_p(t,s)\le e_q(t,s), \qquad \mbox{for all} \quad t\ge s.$$
\end{enumerate}
\end{proposition}

There exist many works where analogs of logarithms were introduced (by several distinct methods), see \cite{AB21}, 
\cite{J08} and \cite{MT09}.

Recall the classical concept of Lyapunov exponents
$$\Upsilon(f)=\limsup\limits_{t\to\infty} \dfrac1t \ln |f(t)|$$
where $f$ is a function (vector function, or matrix function).  

For time scales, instead of using one of the logarithm functions, mentioned above, we can use the following definition \cite{NNHL18}.

\begin{definition} \textit{Lyapunov exponent} of the function $f$ defined on $\mathbb T$ is a real number 
$$\Upsilon_{\mathbb T}(f)=\inf\{a>\nu^*: f(t)/e_a(t)\to 0\}.$$
If the above set is empty, we set $\Upsilon_{\mathbb T}(f)=\nu^*$.
\end{definition}

Even though the Lyapunov exponents depend significantly on the selection of a time scale, the following statement is true.

\begin{proposition}\label{exp0} Let a function $f:{\mathbb R}\to {\mathbb R}$ be such that $\Upsilon(f)\le 0$. Then for any syndetic time scale ${\mathbb T}$ we have $\Upsilon_{\mathbb T}(f)\le 0$. Conversely, if a function $g$ is such that $\Upsilon_{\mathbb T}(g)\le 0$, it may be extended to a function 
${\Hat g}:{\mathbb R}\to {\mathbb R}$ such that ${\Hat g}=\mbox{const}$ on every connected component of ${\mathbb R}\setminus {\mathbb T}$ and 
$\Upsilon_{\mathbb T}(g)\le 0$.
\end{proposition}

This statement is evident.

Given a linear homogeneous system \eqref{linear} with a bounded matrix $A(t)$, and its fundamental system of solutions 
$\Phi=(\varphi_1(t),\ldots,\varphi_n(t))$, recall that $n$ stands for the dimension of the phase space. Let $\Upsilon (\varphi_k)=\alpha_k$, $k=1,\ldots,n$. We say that the fundamental system of solutions is \textit{normal} if it minimizes the sum $\alpha_1+\ldots+\alpha_n$. The existence of a normal fundamental matrix is proved in \cite{NNHL18}, and the proof of that statement is similar to the case of ordinary differential equations \cite{A95}.

Let us recall a classical statement of linear ODE theory. Consider a system on ${\mathbb R}^+$
\begin{equation}\label{odelin}
\dot x=A(t)x.
\end{equation}

Then the so-called Lyapunov inequality takes place \cite{A95}: 
\begin{equation}\label{lyapunov}
\alpha_1+\ldots+\alpha_n\ge \liminf_{t\to \infty} \dfrac{1}{t} \int_{0}^s \mbox{Tr}\, A(s)\, dt
\end{equation}
Here $\mbox{Tr}\, A$ stands for the trace of the matrix $A$. System \eqref{odelin} is called \textit{regular} if inequality \eqref{odelin} turns to equality. In this case, the limit on the right-hand side of \eqref{lyapunov} exists.   

For generic syndetic time scales, a similar estimate was established in \cite{NNHL18}, see the following statement. Let
$$\alpha(t)=\dfrac{\det (I+\mu(t)A(t))}{\mu(t)} \qquad \mbox{if} \quad \mu\neq 0$$
and $\alpha(t)=\mbox{Tr}\, A(t)$ if $\mu(t)=0$.

\begin{theorem} \cite[Theorem 22]{NNHL18} 
For any syndetic time scale $\mathbb T$ and any linear system \eqref{linear}, we have
\begin{equation}\label{lyapt}
\Upsilon[e_{\alpha_1\oplus \alpha_2 \oplus \ldots \oplus \alpha_n}(\cdot,t_0)]\ge \Upsilon[e_{\alpha}(\cdot,t_0)]
\end{equation}
\end{theorem}

It was mentioned in the quoted paper that the question of when the inequality \eqref{lyapt} turns into equality, is open even for systems with constant matrics $A$. 

\begin{definition} The function $f$ is said to have exact Lyapunov exponent (shortly, exact exponent) $\alpha$ if
$$\lim_{t\to\infty} \dfrac{|f(t)|}{e_{\alpha\oplus \varepsilon(t,t_0)}}=0
\qquad \mbox{and} \qquad 
\lim_{t\to\infty} \dfrac{|f(t)|}{e_{\alpha\ominus \varepsilon(t,t_0)}}=\infty.$$
for any $\varepsilon>0$.
\end{definition}

The following statement was proved in \cite{NNHL18} (Theorem 25).

\begin{theorem}\label{texample} If for any eigenvalue $\lambda_i$ of a constant matrix $A$ the function $e_{\lambda_i}(\cdot,t_0)$ has the exact Lyapunov exponent $\alpha_i$ then
\begin{equation}\label{quasiregular}
\Upsilon[e_{\alpha_1\oplus \alpha_2 \oplus \ldots \oplus \alpha_n}(\cdot,t_0)]=\Upsilon[e_{\alpha}(\cdot,t_0)].
\end{equation}
\end{theorem}

\begin{remark} Unfortunately, the systems, satisfying Eq.\, \eqref{quasiregular} can hardly be called regular for the following reason. If  
$$\limsup_{t\to \infty} \mu(t)>0,$$ then adding sufficiently many equations $x_i^\Delta=-x_i$ to any system \eqref{linear}, we get $\nu^*$ in the left and in the right-hand side of Eq.\, \eqref{quasiregular}.
\end{remark}

This is why we introduce another approach to Lyapunov regularity in the next section.

\section{Regularity of finite-dimensional systems on time scales}

In this section, we always deal with time scales such that $t_+=+\infty$. Without loss of generality, we assume that $t_-=0$. Besides, we assume that $\dim X<+\infty$. Given a fundamental matrix $\Phi=(\varphi_1,\ldots,\varphi_n)$ of system \eqref{odelin}, we consider the value 
$$S(\Phi)=\sum_{i=1}^n \Upsilon (\varphi_i).$$
We define 
$$S(A)=\min_\Phi S(\Phi).$$ 
In this classical case, system \eqref{odelin} defined on $[0,\infty)$ is called \emph{regular} if 
$$S(A)-\liminf_{t\to+\infty} \int_{0}^t \mathrm{Tr}\, A(s)\, ds=0$$
(the left-hand side of the latter formula is always non-negative).

The following statement was proved by Basov \cite{B52}. 

\begin{theorem}\label{th_byl1} System \eqref{odelin} of ordinary differential equations is regular if and only if there exists a constant matrix $$B=\mathrm{diag}\, (b_1,\ldots,b_n)$$ and a continuous non-degenerate matrix $L(t)$ with $\Upsilon L=\Upsilon {L^{-1}}=0$ such that the substitution $x=L(t)y$ reduces system \eqref{odelin} to system $\dot y=By$.
\end{theorem}

We use the latter statement to define regular systems on time scales. 

\begin{definition}\label{defbyl2} Let the phase space $X$ be a Euclidean space ${\mathbb R}^n$. System \eqref{linear} is regular if and only if there exists a constant diagonal matrix/operator 
$$B=\mathrm{diag}\, (b_1,\ldots,b_n)$$ 
and a piece-wise continuous non-degenerate transformation $L(t)$ with 
$$\Upsilon L=\Upsilon {L^{-1}}=0$$ 
such that the substitution $x=L(t)y$ reduces system \eqref{linear} to system 
\begin{equation}\label{yb}
y^\Delta=\Hat B(t)y. 
\end{equation}
Here 
\begin{equation}\label{barb}
\begin{array}{c}
\Hat B(t)=B \quad \mbox{if} \quad \mu(t)=0 \quad \mbox{and}\\
\Hat B(t)=\dfrac{\exp(B\mu(t))-E}{\mu(t)}\quad \mbox{if} \quad \mu(t)>0.
\end{array}
\end{equation}
\end{definition}

Observe that system \eqref{yb} has the Cauchy matrix $\Phi(t,\tau)=\exp (B(t-\tau))$. 
Applying methods of the Floquet theory, one can easily obtain that periodic systems on periodic time scales with periods coinciding are regular. Besides, the system with a constant matrix is regular provided all the solutions have exact exponents. This is true even for non-periodic time scales (see Theorem \ref{texample}). The latter fact can be proved similarly to the quoted result of \cite{B52}.

Now we study the existence of bounded solutions for systems with regular linear parts. 

\begin{theorem}\label{th_r_1} Let system \eqref{linear} be regular on ${\mathbb R}^+\bigcap {\mathbb T}$ where the time scale ${\mathbb T}$ is syndetic. Then for any $\alpha, \beta, \gamma>0$ there exist positive constants $c_{1,2}>0$ and $h>0$ such that if
\begin{equation}\label{eqa}
\|a(t,x)\| \le c_1 \|x\|^{1+\alpha} +c_2 \exp(-\beta t) \|x\|+h \exp(-\gamma t),
\end{equation}
then system \eqref{timescale1} has a bounded solution that is exponentially decaying at infinity.
\end{theorem}

\noindent\textbf{Proof.} 

First of all, we consider a set of linear non-homogeneous systems 
\begin{equation}\label{nonhom}
x^\Delta=A(t)x+f(t)
\end{equation}
where 
\begin{equation}\label{clambda}
|f(t)|\le C\exp (-\lambda t).
\end{equation}
Let ${\mathcal C}_\lambda$ be the space of functions that satisfy Eq.\eqref{clambda} for some $\lambda$. In this space, we introduce the norm 
$$\|f\|_\lambda=\|f\|_{\lambda,{\mathbb T}}=\sup_{t\in {\mathbb T}} |f(t)|\exp(\lambda t).$$

Select a $\lambda\in (\max(\gamma/(1+\alpha),\gamma-\beta),\gamma)$ where $\alpha$, $\beta$ and $\gamma$ are defined by formula \eqref{eqa}. Then $\lambda (1+\alpha)>\gamma$ and $\lambda+\beta>\gamma$.

\begin{lemma} Let a system \eqref{linear} be regular on a syndetic time scale $\mathbb T$. Then for any $\gamma>\lambda$ there exists a constant $C_{\gamma,\lambda}>0$ and a linear operator ${\mathcal L}_\gamma: {\mathcal C}_\gamma \to {\mathcal C}_\lambda$ mapping the right-hand side $f$ of system \eqref{nonhom} to its bounded solution with $\|{\mathcal L}_\gamma\|\le C_{\gamma,\lambda}$.
\end{lemma}

\noindent\textbf{Remark}. For non-syndetic time scales, the above statement is, generally speaking, wrong. For instance, the equation $x^{\Delta}=e^{-t}$ does not have a bounded solution on the time scale $\{3^{3^n}: n\in {\mathbb N}\}$.

\noindent\textbf{Proof of Lemma}. Evidently, it suffices to prove the result for a diagonal system with coefficients $\Hat b_j(t)$, being defined by a constant $b_j$ by formulae similar to \eqref{barb}. Moreover, one can only consider scalar equations of the type $u^\Delta = {\Hat b}(t)u+f(t)$ where the homogeneous system $u^\Delta = {\Hat b}(t) u$ has a solution $\exp(bt)$.

Given an rd-continuous function $f$ we define the extension ${\Hat f}:{\mathbb R}\mapsto X$ as follows: ${\Hat f}(t)=f(t)$ if $t\in {\mathbb T}$, $\mu(t=0)$ and 
\begin{equation}\label{ebarf}
{\Hat f}(t)=\mu(\Hat t) b f(\Hat t)/(1-\exp(-b\mu(\Hat t)))
\end{equation}
otherwise. Here $\Hat t=\sup\{s\in {\mathbb T}: s\le t\})$. Equation \eqref{ebarf} is an equivalent form of the equality
$$\int_{{\Hat t}}^{\sigma({\Hat t})} e^{-bs} \Hat f \, ds= 
\int_{{\Hat t}}^{\sigma({\Hat t})} e^{-bs} f(s) \, \Delta s.$$
Since the considered time scale $\mathbb T$ is syndetic, there exist a constant $K$ that depends on the matrix $A$ only and is such that for any $b\in \{b_1, \ldots, b_n\}$ and any $f$ we have $\|{\Hat f}\|_{\lambda,{\mathbb R}^+} \le 
K|f|_{\lambda,{\mathbb T}}$. 

Then the required solution may be found by the formula 
${\mathcal L}f=h$ where 
$$h(t)=\int_{0}^t \exp(b(t-s)) f(s)\, \Delta s=\int_0^t \exp(b(t-s)) {\Hat f}(s)\, ds$$
for any $t\in {\mathbb T}$ if $\lambda\ge -b$.

Similarly,
$$h(t)=-\int_t^\infty \exp(b(t-s)) {\Hat f}(s)\, ds$$
otherwise.
$\square$

\medskip

Now, we start proving Theorem \ref{th_r_1}. Observe that having inequality \eqref{eqa}, we can replace any of constants $\alpha$, $\beta$, or $\gamma$ with a smaller value, and the inequality will still be true. Suppose, without loss of generality, that the segment $(-\gamma,0)$ does not contain any exponents of system \eqref{odelin} (the total number of these exponents cannot exceed the dimension of the system).

We demonstrate that the considered system has a solution from ${\mathcal C}_\lambda$ where $\lambda$ was selected above.

First of all, observe that, given a function $x(\cdot)\in {\mathcal C}_\lambda$, 
$$\|a(\cdot, x(\cdot))\|_\gamma \le (c_1+c_2)\|x(\cdot)\|_\lambda+h$$
and 
$$\|{\mathcal L}_\gamma a(\cdot, x(\cdot))\|_\lambda \le C_{\gamma,\lambda}(c_1+c_2)\|x(\cdot)\|_\lambda+C_{\gamma,\lambda}h.$$

Secondly, the operator ${\mathcal L}_\gamma$ acting from ${\mathcal C}_\gamma$ to ${\mathcal C}_\lambda$ is compact. Indeed, let $B$ be a unit ball in  ${\mathcal C}_\gamma$. Then norms of functions ${\mathcal L}_\gamma f$ where $f\in B$ are uniformly bounded in ${\mathcal C}_\lambda$ as well as their derivatives. Moreover, the set ${\mathcal L}_\gamma B$ is convex.

Let $\kappa=C_{\gamma,\lambda}h/(C_{\gamma,\lambda}(c_1+c_2))$. All in all, the nonlinear operator 
${\mathcal L}a(t,\cdot)$ maps a $\kappa$ ball in ${\mathcal C}_\lambda$ centered at zero to it compact subset. Then, applying Schauder's fixed point theorem, we obtain the existence of the desired solution. 
$\square$

\section{\bf Conclusion and discussion}

In this paper, we studied the existence of a bounded solution of an 'almost linear' time scale dynamics system. First of all, we consider the condition of exponential dichotomy (hyperbolicity) on a family of subsegments of the time scale. Regarding this condition, we formulate and prove the condition sufficient for the existence of a bounded solution. Later on, we study the general condition for the solvability of an 'almost linear' time-scale system. As a particular case, we discuss uniformly hyperbolic and regular systems on time scales.

The following two remarks indicate the possible directions for future work.

\noindent\textbf{Remark.} The statement of Theorem \ref{th_r_1} admits various generalizations that are supposed to be done in the future.
\begin{enumerate}
\item Similarly to Theorem \ref{th_hyp_2}, one can consider systems, regular on both half-lines, positive and negative. In this case, the nonlinearities must be exponentially decaying in both directions.
\item Bounded solutions of systems, consisting of regular and hyperbolic blocks, can also be studied.
\item 'Exponentially small' perturbations of time scales can also be considered.
\end{enumerate}

\noindent\textbf{Remark}. The case of the so-called $\nabla$ derivative can also be considered. In this case, a forward non-uniqueness of solutions may be there. In fact, the system with $\nabla$ -- derivatives can be reduced to a system with a $\Delta$ -- derivative by applying the transformation $t\mapsto -t$.

\section*{Acknowledgements}
The work of the second co-author was supported by Gda\'{n}sk University of Technology by the DEC 14/2021/IDUB/I.1 grant under the Nobelium - 'Excellence Initiative - Research University' program. Coauthors are grateful to Prof. S.\,Yu.\, Pilyugin for fruitful discussion and useful advice.


\begin{thebibliography}{99}
\bibitem{A95} L.\, Ya.\, Adrianova, \emph{Introduction to Linear Systems of Differential Equations}, American Mathematical Society, 1995. 
\bibitem{AB21} D.\, R.\, Anderson, M.\, Bohner, \emph{A multivalued Logarithm on Time Scales}, Applied Mathematics and Computation, \textbf{397}, 15 May 2021, 125954.
\bibitem{AH90} B.\, Aulbach, S.\, Hilger, \emph{Linear Dynamic Processes with Inhomogeneous Time Scale}, In Nonlinear Dynamics and Quantum Dynamical Systems (Gaussig, 1990), volume 59 of Math. Res., pages 9--20. Akademie Verlag, Berlin, 1990.
\bibitem{B52} V.\, P.\, Basov, \emph{On the structure of a solution of a regular system}, Vestnik Leningr. Univ. (1952), no. 12, 3--8. (Russian)
\bibitem{BP01}  M.\, Bohner, A.\, Peterson, \emph{Dynamic Equations on Time Scales: An Introduction with Applications}, Birkh\"auser, Boston, 2001.
\bibitem{BP03} M.\,Bohner, A.\,Peterson, \emph{Advances in Dynamic Equations on Time Scales}. Birkh\"{a}user Boston Inc., Boston, MA, 2003.
\bibitem{C12}  J.\, L.\, Cie\'sli\'nski, \emph{New Definitions of Exponential, Hyperbolic, and Trigonometric Functions on time scales}, J. Math. Anal. Appl. \textbf{388} (2012) 8--22.
\bibitem{C78} W.\, A.\, Coppel, \emph{Dichotomies in Stability Theory}, Lecture Notes in Mathematics No. 629, Springer-Verlag, Berlin, 1978.
\bibitem{G19} S.\, Georgiev, \emph{Functional Dynamic Equations on Time Scales}, Springer-Verlag, 2019, ISBN 978-3-030-15419-6.
\bibitem{J08} Jackson, B., \emph{The Time Scale Logarithm}, Appl. Math. Letters \textbf{21}, (2008) 215--221.
\bibitem{KV06} S.\, Kryzhevich, V.\, Volpert, \emph{On Different Types of Solvability Conditions for Differential operators}, Electronic Journal of Differential Equations, 2006, vol.2006, 1--24.
\bibitem{M16} A.\,A\, Martynyuk, \emph{Stability Theory for Dynamic Equations on Time Scaless, Systems
and Control}:, Birkh\"auser/Springer, [Cham], 2016.
\bibitem{MT09} D.\, Mozyrska, D.\,F.\,M\, Torres, \textit{The Natural Logarithm on Time Scales}, J. Dyn. Sys. Geom. Theories 7:1 (2009) 41--48.
\bibitem{NNHL18} K.\,C.\,Nguen, T.\,V.\, Nhung, T.\,T.\, Anh Hoa, N.\,C.\,Liem, \emph{Lyapunov Exponents for Dynamic Equations on Time Scales}, Dynamic Systems and Applications, \textbf{27}, No. 2 (2018), 367--386.
\end{thebibliography}
\end{document}